\documentclass{birkmult}
 \newtheorem{thm}{Theorem}[section]

 \theoremstyle{definition}
 
 \theoremstyle{remark}
 \newtheorem{rem}[thm]{Remark}
 
 \numberwithin{equation}{section}

\usepackage{amssymb,amsmath,amsfonts}

\def\G{\mathcal{G}}

\def\L{\mathcal{L}}
\def\P{\mathcal{P}}
\def\T{\mathcal{T}}

\begin{document}
%
%
%
%
%
%
%
%
%
\title[Hydrodynamics, probability and the geometry of the diffeomorphisms group]
 {Hydrodynamics, probability and the geometry of the diffeomorphisms group}
\author[Ana Bela Cruzeiro]{Ana Bela Cruzeiro}

\address{%
Dep. Matem\'atica I.S.T. (TUL)\\
Av. Rovisco Pais\\
1049-001 Lisboa\\
Portugal\\
\&\\
Grupo de F\'{\i}sica-Matem\'atica da Universidade de Lisboa\\
Av. Prof. Gama Pinto, 2\\
1649-003 Lisboa\\
Portugal}

\email{abcruz@math.ist.utl.pt}

\thanks{This work was completed with the support of  the \textit{Funda\c{c}\~ao para a Ci\^encia e a Tecnologia} through the project PTDC/MAT/69635/2006.}

\subjclass{Primary 37L55; Secondary 35Q30, 58E30, 58J65, 60J60, 76D05}

\keywords{Navier-Stokes equation; diffeomorphisms group; stochastic geodesic equation.}

\date{January 1, 2004}

\begin{abstract}
We characterize the solution of Navier-Stokes equation as a stochastic  
geodesic on the diffeomorphisms group, thus generalizing Arnold's  
description of the Euler flow.
\end{abstract}

\maketitle
\section{Introduction}

Laplacian canonical determinism (quoting Laplace himself \textit{``We may regard the present state of the universe as the effect of its past and the cause of its future..."}) has been called into question by the studies of H. Poincar\'e and  J. Hadamard of complex highly unstable trajectories of a system of three bodies moving under the effect of gravitation forces. Much later E. Lorenz considered some simplified model for atmospheric convection with a small number of degrees of freedom and also found
\textit{chaotic} motions, the term chaotic meaning that small changes in the initial conditions  produce
 big perturbations in the  behaviour of the trajectories of motion (sometimes popularly referred to as the ``butterfly effect").

\indent Dynamics may therefore exibit a  chaotic (random?) behaviour even 
when the governing equations  are deterministic. Beyond classical dynamics, the most convincing illustration
of this is, of course, quantum mechanics.

\indent In Turbulence, modelized by the (deterministic) incompressible Navier-Stokes equation 
\begin{equation*}
\frac{\partial u}{\partial t} +u.\nabla u =\nu \Delta u -\nabla p,\;div~ u=0,\;u(0,x)=u_0 (x)
\end{equation*}
with small viscosity coefficient $\nu >0$, such classical chaotic phenomena really occur. On the other hand statistical effects also manifest themselves: they are called \textit{coherent structures} (for example the Jupiter's red spot
 or A\c{c}ores anticyclone). We refer to \cite{R} for an interesting discussion of these problems.
 
\indent There are various possible reasons to look for probabilistic content  in Fluid Dynamics. The obvious observations that initial conditions are difficult to know precisely and  that the observable computed in practice have to be  statistical quantities points out to some statistical approach, found, for example, in works of Vishik, Fursikov, Foias, etc. (\cite{F}).

\indent One can of course add  some noise to given deterministic systems and study the randomly perturbed ones: then  the resulting properties depend on the chosen perturbation and these studies do not describe the chaotic behaviour of the deterministic equations themselves. These approaches have known an impressive development in the recent years.
   
\indent The intrinsic instability of the hydrodynamical motions, more precisely  those modeled by the Euler equations (the case where the viscosity is zero) was explained geometrically without using  Probability theory by V. Arnold (\cite{A}). In his famous paper Arnold  showed that the solution of Euler equation 
\begin{equation*}
\frac{\partial u}{\partial t} +u.\nabla u =-\nabla p, \; div\, u=0
\end{equation*}
corresponds to the velocity of a flow which is critical for the action functional
\begin{equation*}
S[g]=\frac{1}{2} \int ||\dot g (t)||^2_{L^2}dt
\end{equation*}
where $g$ denote volume preserving diffeomorphisms of the underlying manifold. Accordingly
Euler equation can be regarded as the   geodesic equation for the $L^2$ metric and one derives the instability of the motion strictly from geometrical arguments, namely from the negativity of some curvature of the space of diffeomorphisms. This is not in contradiction with the formation of the abovementionned coherent structures since Arnold describes  the so-called \textit{Lagrangian} flow, namely $g_t (x)$ such that 
 \begin{equation*}
 \partial_t g_t (x)= u (t, g_t (x))
 \end{equation*}
and not the usual \textit{Eulerian} flow (of velocities). It seems that there is a kind of chaotic behaviour of the Lagrangian flow and a more predictable character of the Eulerian flow $u(t,.)$.

\indent We  note that this ``principle" seems to apply as well to properties like ergodicity. Indeed,
there are reasons to believe that the flow of velocities is ergodic (cf. \cite{B} or a more recent discussion in \cite{G}), invariant measures can be constructed (cf. \cite{A-C} and references therein for the two-dimensional periodic case). On the other hand the Lagrangian flow is not ergodic; this has been proved in \cite{C-M}, for the periodic $n$-dimensional case.

\indent The statistical effects observed in the evolution of the system are, on the other hand,
 related to the large number of degrees of freedom in presence, a situation quite distinct
 from the dynamical systems studied by Poincar\'e or Lorenz.

\indent Following an idea of Yasue etal (\cite{N,Y}), consisting in interpreting the Laplacian term of Navier-Stokes equation as  origin of the random behaviour of fluid particles and the drift as an avereged velocity on the stochastic (Lagrangian) motion, one is able to 
\begin{itemize} 
\item[{\rm a)}] obtain an intrinsic probabilistic representation of  Navier-Stokes solution;
\item[{\rm b)}] derive it as a solution of a (stochastic) variational principle;
\item[{\rm c)}] deduce  the corresponding (stochastic) equation for geodesics, the associated Jacobi equation and therefore  consider the  stability properties of the motion.
 \end{itemize}

\indent The methodology thus generalize probabilistically Arnold's approach  to the case where viscosity is not zero.

\section{Geometry of the diffeomorphisms group}

Let $G_V$ denote the  group of volume preserving homeomorphisms on a manifold. Throughout this paper we consider this manifold to be $T^2$, the two dimensional torus, for simplicity, since we can as well treat the torus in any dimension (cf. the remark after Theorem \ref{thm2}). On a  general manifold the formulation of our results may be a more delicate matter.

\indent The tangent space of the group $G_V$ will be denoted by $\G_V$; it consists of vector fields with zero divergence. The $L^2$ norm defines on 
$\G_V$ a canonical Hilbert structure. More generally, if $G^s = \{g: T^2 \rightarrow T^2 ~ \hbox{bijection},~ g, g^{-1} \in H^s\}$ then for  $s>2$, ~ $G^s$ is a $C^{\infty}$ infinite dimensional Hilbert manifold (cf.  \cite{E-M, S}). The associated tangent space, endowed with the $H^s$ norm, will be denoted by $\G^s$ and the tangent space for the corresponding volume preserving subgroup $G_V^s$ by $\G_V^s$ .

\indent We shall consider Fourier developments, namely, for $f\in L^2$, 
\begin{equation*}
f(\theta)=2 \sum_{k\in {\tilde Z}^2}~\Re{\hat f}(k)~\cos( k.\theta)-\Im{\hat f}(k)~\sin( k.\theta) 
\end{equation*}
where $\hat f$ is the Fourier transform of $f$,
\begin{equation*}
{\hat f}(k) = \frac{1}{(2\pi )^2} \int_{T^2} ~exp ~i(k.~\theta) d\theta
\end{equation*}
and
\begin{equation*}
\tilde Z^2 = \{ k= (k_1 ,k_2 )\in Z^2 : k_1 >0 ~ \hbox{or} ~  k_1 =0, k_2 >0\}
\end{equation*}
Notice that $f$ is real iff ${\hat f}(-k)=\bar{\hat f}(k)$.

\indent We take for orthonormal  basis of  $\G_V^s$  the following vector fields:
 \begin{align*}
 A^s_k & =\frac{1}{|k|^{s+1}}[(k_2\cos k.\theta) \partial_1-(k_1\cos k.\theta)\partial_2)]\\
 B^s_k &= \frac{1}{|k|^{s+1}}[(k_2\sin k.\theta) \partial_1-(k_1\sin k.\theta)\partial_2)]
 \end{align*}
 $k\in {\tilde Z^2}-\{(0,0)\}$, where $|k|^2=k_1^2+k_2^2$ and $\partial_i={\frac{\partial}{\partial \theta^i}}$, together with the constant vector fields.

When $s=0$ we  omit the corresponding index in the notation.

\indent The  constants of structure  of $\G_V^s$ are
\begin{align*}
[A^s_k ,A^s_l ]&=\frac{[k,l]}{2|k|^{s+1}|l|^{s+1}}(|k+l|^{s+1}B^s_{k+l}+|k-l|^{s+1}B^s_{k-l}) \\
[B^s_k ,B^s_l ]&=-\frac{[k,l]}{2|k|^{s+1}|l|^{s+1}}(|k+l|^{s+1}B^s_{k+l}-|k-l|^{s+1}B^s_{k-l})\\
[A^s_k ,B^s_l ]&=-\frac{[k,l]}{2|k|^{s+1}|l|^{s+1}}(|k+l|^{s+1}A^s_{k+l}-|k-l|^{s+1}A^s_{k-l})\\
 [\partial_i ,A^s_k]& =-k_i B^s_k\\
 [\partial_i ,B^s_k]& =k_i A^s_k
  \end{align*}

\indent One can find a detailed proof of these expressions in \cite{C-F-M} for the case where $s=0$, a computation that extends easily to the general case. For the torus of higher dimensions we refer to \cite{C-M}.

\indent Define  the following functions
\begin{align*}
\alpha^s_{k,l}: & = \frac{1}{4(|k||l||k+l|)^{s+1}}(|k+l|^{2(s+1)}-|k|^{2(s+1)}+|l|^{2(s+1)}) \\
\beta^s_{k,l}: & = \frac{1}{4(|k||l||k-l|)^{s+1}}(|k-l|^{2(s+1)}-|k|^{2(s+1)}+|l|^{2(s+1)}) \\
[k,l] & = k_1 l_2-k_2 l_1
\end{align*}

\indent The Christoffel symbols for the Levi-Civita connection associated to the $H^s$ metric can be explicitly computed; they are given by
\begin{align*}
\Gamma^s_{A^s_k ,A^s_l}&=[k,l]( \alpha^s_{k,l} B^s_{k+l}+\beta^s_{k,l} B^s_{k-l})\\
\Gamma^s_{B^s_k ,B^s_l}&=[k,l]( -\alpha^s_{k,l} B^s_{k+l}+\beta^s_{k,l} B^s_{k-l})\\
\Gamma^s_{A^s_k ,B^s_l}&=[k,l]( -\alpha^s_{k,l} A^s_{k+l}+\beta^s_{k,l} A^s_{k-l})\\
\Gamma^s_{B^s_k ,A^s_l}&=[k,l]( -\alpha^s_{k,l} A^s_{k+l}-\beta^s_{k,l} A^s_{k-l})
\end{align*}

\indent In the case where $s=0$ the proof can be found in \cite{C-F-M}.

\indent A similar computation for the space $Diff (S^1 )$ was done in \cite{A-M}.

\indent In order to study the geodesics corresponding to the Euler equation 
in fact only this case is (physically) relevant: the $L^2$ norm (i.e., the energy) has to be considered, even if one is lead,
 for technical reasons, to consider the equations in higher order Sobolev spaces, where 
 the $L^2$ norm gives rise to weak Riemannian structures (cf. \cite{E-M},\cite{S}).

\indent Concerning the corresponding Ricci curvature (for the $H^s$ metric) we have the
  following explicit expression   
\begin{align*}
 Ricci_s (A^s_j) =& -2 \sum_i [i,j]^2 (\alpha^s_{i,j} \alpha^s_{j,i} + \beta^s_{i,j} \beta^s_{j,i}) A^s_j  \\
&-2\sum_i [i,j]^2 \frac{|i+j|^{s+1}}{(|i| |j|)^{s+1}}\beta^s_{i+j,i}A^s_j -2\sum_i [i,j]^2 \frac{|i-j|^{s+1}}{(|i| |j|)^{s+1}}\beta^s_{i-j,i}A^s_j 
\end{align*}
and an analogous one for $Ricci_s (B^s_j)$ (the curvature is a diagonal operator).
 
\indent This result was proved in \cite{C-M} for the case $s=0$. We notice that, in this case,
the last two terms in each of the Ricci curvature expressions are symmetric,
\begin{equation*}
\frac{|i+j|}{|i| |j|} \beta_{i+j,i}+\frac{|i-j|}{|i| |j|} \beta_{i-j,i} =\frac{-(i|j)}{|i|^2 |j|^2} +\frac{(i|j)}{|i|^2 |j|^2}=0
\end{equation*}
and we are there left with only the first terms, which reduce to:
\begin{align*}
Ricci (A_j )&= -\sum_i [i,j]^4 \frac{(|i|^2 +|j|^2 )}{|i|^2 |j|^2 |i-j|^2 |i+j|^2} A_j  \\
Ricci (B_j )&= -\sum_i [i,j]^4 \frac{(|i|^2 +|j|^2 )}{|i|^2 |j|^2 |i-j|^2 |i+j|^2} B_j.
\end{align*}

\indent We observe that the $Ricci_s$ tensors are given by divergent quantities.

\indent On the Lie algebra $\G^s_V$ we may define the following process
\begin{equation*}
d x(t)=\sum_{k\neq 0}  (A^s_k ~dx_k (t) +B^s_k ~d y_k (t))
\end{equation*}
with $s\geq 2$, where $x_k, ~y_k$ are independent copies of real Brownian motions. The stochastic flow
\begin{align*}
d g(t)&=(od x(t))(g(t))\\
g(0)&=e
\end{align*}
($e$ being the identity)  is a well defined stochastic flow of diffeomorphisms  (cf. \cite{K} and  \cite{M, C-C} for exponents $s$ non necessarily greater then $2$). Its generator is
\begin{equation*}
L^s =\frac{1}{2}\sum_{k\neq 0} \partial_{A^s_k}\partial_{A^s_k} F(g)+\frac{1}{2}\sum_{k\neq 0}\partial_{B^s_k}\partial_{B^s_k} F(g).
 \end{equation*}

\begin{thm}
The process $g(t)$ coincides with the Brownian motion associated to the metric $H^s$.
\end{thm}

\indent The proof follows from the explicit expression of the Christoffel symbols, more precisely from the fact that $\Gamma^s_ {A^s_k , A^s_k}=\Gamma^s_ {B^s_k , B^s_k}=0$.

\vskip 5mm

\indent The Laplace-Beltrami
 operator on vector fields $Z$ reads,
\begin{equation*}
\L^s Z =\frac{1}{2}\sum_{k\neq 0}  \nabla^s_{A^s_k} \nabla^s_{A^s_k} Z +
\frac{1}{2} \sum_{k\neq 0} \nabla^s_{B^s_k} \nabla^s_{B^s_k} Z. 
\end{equation*}
where $\nabla^s$ denotes the covariant derivative for the $H^s$ metric.
 
\indent  Then we have the following fundamental result that enables us to relate the classical Laplacian term
 in Navier-Stokes equations to the stochastic geometry on the space of (volume preserving)
 diffeomorphisms:

\begin{thm}{\rm(\cite{C-C})}\label{thm2}
When computed on functionals of the form $F(g) (\theta )=f (g(\theta ))$, $f$ smooth, we
have
\begin{equation*}
L^s F (g)= \frac{c}{2} \Delta f 
\end{equation*}
with $c=\sum_{k\neq 0}  \frac{(k_i )^2}{|k|^{2 (s+1)}}$.
\end{thm}

\indent The proof can be found in \cite{C-C} (c.f. also \cite{C-S}).

\begin{rem}
We  mention  that most of the results stated in this paper extend without too much trouble (but with a certain amount of work) to the case of a torus
of dimension higher then two, once we choose a suitable basis of $\G_V$ and the corresponding ones for the higher order Sobolev spaces.

\indent Denote by $\tilde Z^d$ a subset of $Z^d$ such that each equivalence class of the equivalence relation defined by $k\simeq k'$ if $k+k' =0$ has a unique representative in $\tilde Z^d$. 
We should take, for every $k\neq 0$, an orthonormal basis of $\{x\in R^d: k.x=0 \}$, namely 
$\epsilon_{k,\alpha},~\alpha=1,...,d$ and define
\begin{align*}
A_{k,\alpha } & =\epsilon_{k,\alpha } cos (k.\theta ) \\
B_{k,\alpha } & =\epsilon_{k,\alpha } sin (k.\theta )
\end{align*}
to obtain (after normalization) the corresponding basis of $\G^s_V$.

\indent Furthermore we should assume that $\sum_{k, \alpha } [(\epsilon_{k,\alpha} )^i ]^2$ is independent of $i$ and that $\sum_{k, \alpha } [(\epsilon_{k,\alpha} )^i )(\epsilon_{k,\alpha} )^j )=0$ for $i\neq j$ in order to obtain the generalization of last theorem.

\indent For the generalization of some of the results in this paper to torus of dimensions bigger then $2$ we refer to \cite{C-M} and to \cite{C-S}.
\end{rem}

\indent For regularity reasons, in order to avoid geometrical divergent quantities such as the $Ricci$ tensors, we are going from now on to consider a truncated Brownian motion $x^N (t)$ defined as

\begin{equation*}
d x^N(t)=\sum_{0<|k| \leq N }  (A^s_k ~dx_k (t) +B^s_k ~d y_k (t))
\end{equation*}
and the corresponding process $g^N (t)$ with generator $L^{s,N}$.
Then, by the same arguments of Theorem 2.2. we can show that,
when computed on functionals of the form $F(g) (\theta )=f (g(\theta ))$, $f$ smooth, we
have
\begin{equation*}
L^{s,N} F (g)= \frac{c_N }{2} \Delta f 
\end{equation*}
with $c_N=\sum_{0<|k|\leq N}  \frac{(k_i )^2}{|k|^{2 (s+1)}}$.

\section{A variational principle}

For a time dependent vector field $u\in L^2 ([0,T]; \G_V)$, let $g^N_u$ be a solution of the s.d.e. with values in $G_V$,

\begin{equation*}
dg^N_u (t) =(u(t) dt +\sqrt\frac{\nu}{c_N} od x^N(t))g^N_u (t), ~~~ g_u (0)=e
\end{equation*}

with generator 
\begin{equation*}
L_u^{s,N} F (g)=\frac{\nu}{c_N} L^{s,N} F +\sum_k u^{A^s_k} \partial_{A^s_k} F+\sum_k u^{B^s_k}
\partial_{B^s_k}F.
\end{equation*}

\indent A weak solution of such stochastic differential equation can be shown to exist (cf. \cite{C-C}).

\indent Then, when computed on functionals of the form $F(g) (\theta )=f (g(\theta ))$, $f$ smooth, we have
\begin{equation*}
L^s_u F (g)= \nu \Delta f +u.\nabla f.
\end{equation*}
Notice that the constant appearing in the diffusion coefficient is chosen to produce, after ``projection" on pointwise functions,  the viscosity parameter $\nu$ in front of the Laplacian.

\indent For smooth functionals $F(g,t)$  denote
\begin{equation*}
D F(g(t),t)=lim_{\epsilon \rightarrow 0} E^{\P_t} (F(g(t +\epsilon ), t+\epsilon
)-F(g(t ),t))
\end{equation*}

\indent For (time dependent) vector fields $Z(t,.)\in \G_V$, and
using the It\^o
parallel transport (associated to the Levi-Civita connection) $\T_{s\leftarrow t}$ above the process $g(t)$, which is the solution of
 the stochastic differential equation
\begin{align*}
d_ t [\T_{s\leftarrow t} ] Z_t &=- (\Gamma_k )_{g(t)}~o~d x_t^k ~[\T_{s\leftarrow t } ]Z_t \\
[\T_{s\leftarrow s} ]  &= Id
\end{align*}
define
\begin{equation*}
D Z (g(t),t)=lim_{\epsilon \rightarrow 0} \frac{1}{\epsilon} E^{\P_t} ( \T_{t \leftarrow
(t+\epsilon )} Z(g(t+\epsilon ) , t+\epsilon  )- Z(g(t),t))).
\end{equation*}

\indent Consider  now the following (stochastic) action functional
\begin{equation*}
S[g]=\frac{1}{2} E\int_0^T ||Dg (t)||^2_{L^2} dt -\frac{1}{2} E||Dg(0)||^2_{L^2}
\end{equation*}

\indent The Navier-Stokes solution can be characterized as the drift of a critical process for this action. More precisely we have,

\begin{thm} {\rm(\cite{C-C})}

The process $g^N_u (t)$ is a critical point of the energy functional $S$ if and only if $u$ verifies the N. S. equation
\begin{align*}
\frac{\partial u}{\partial t} +u.\nabla u =\nu \Delta u -\nabla p\\
div~ u=0, ~~ u(T, \theta )=u_T (\theta )
\end{align*}

\end{thm}

\indent We precise that \textit{critical} here means  in the space  of continuous $G_V$-valued semi-martingales $g(t)$ with $g(0)=e$ with variations defined as
\begin{equation*}
(D_l )S [e_. (v)] =\frac{d}{d \epsilon } S[e_. ( v) o g(. )]|_{\epsilon =0}
\end{equation*}
where
\begin{align*}
e_t (v), v \in C^1 ([0,T]; \G_V^{\infty} ), ~v(0, .)=0 \\
\frac{de_t (v)}{dt}= \dot v (t, e_t (v)),~ e_0 (v)=e.
\end{align*}

\section{Stochastic geodesics equation}

To a variational principle a Euler-Lagrange dynamics can be associated. In this framework the corresponding equations are stochastic and correspond to a ``regularized" version of the equation of geodesics  for the $L^2$ metric.

\indent Writting $e_k$ for a generic element of the basis of $\G_V^s$ ($A_k^{s}$ or $B_k^{s}$) and $u^k =u^{e_k}$, define the following truncated approximation of the $Ricci_s$ tensor

\begin{equation*}
R_s^N (\tilde{u}) =\sum_{0<|k|\leq N} \sum_{l,m}\Gamma^s_{e_k ,e_l}   
(\Gamma^s_{e_m,e_k}|e_l ) \tilde{u}^m -([e_k ,e_m ]|e_l )\Gamma^s_{l,k} \tilde{u}^m
\end{equation*}
\begin{align*}
~~~~~=& -2 \sum_{0<|k|\leq N}\sum_m  [k,m]^2 (\alpha^s_{k,m} \alpha^s_{m,k} + \beta^s_{k,m} \beta^s_{m,k}) \tilde{u}^m  \\
&-4\sum_{0<|k|\leq N} [k,m]^2 \frac{|k+m|^{s+1}}{(|k| |m|)^{s+1}}\beta^s_{k+m,k}{\tilde u}^m
\end{align*}

\begin{thm}
Let $s\geq 0$. The vector field  $\tilde u(t,\theta)=-u((T-t), \theta )$ is a solution of the Navier-Stokes equation if and only if $u$ is a martingale of the process $g^N_u$.

\indent Equivalently,
\begin{equation*}
DD g^N_u (t)-\frac{\nu}{2c} R^N_s (D g^N_u )=0 ~~a.e.
\end{equation*}
with
\begin{equation*}
g^N_u (0)=e,~~~Dg^N_u (T)=u(T,g^N_u (T))
\end{equation*}
or
\begin{equation*}
\frac{\partial \tilde u}{\partial t} +\tilde u^k \nabla_{e_k} \tilde u -\L^{s,N} \tilde{u} +\frac{\nu}{2c_N}
R^N_s (\tilde u)=0
\end{equation*}
\end{thm}

\begin{rem}
\indent
\begin{enumerate}
\item The presence of the truncated $Ricci$ term is due to the choice of the It\^o stochastic parallel transport 
(associated with the Levi-Civita connection). This term disappears if we define the operator $D$ using the so-called ``damped" parallel transport, which incorporates the effect of the curvature.
 \item The covariant derivation appearing in the last formula of the Theorem is the one associated to the $L^2 =H^0$ metric, as it should be, since it is the (Euler) non linear term of the equation. On the other hand the stochastic perturbations of the Euler equation, responsible for the second order term of the Navier-Stokes equation,
 are related to the $H^s$ metric.
 \item The expression
 
\begin{equation*}
 \L^{s,N} \tilde{u} -\frac{\nu}{2c_N}
R^N_s (\tilde u)=\nu \Delta \tilde{u}
\end{equation*}
may be well defined with each of its separate terms divergent, depending on the regularity of the Navier-Stokes velocity function $u$.

\end{enumerate}

\end{rem}

\begin{proof}
The term $\tilde u^k \nabla_{e_k} \tilde u$ corresponds, as said before, to the non-linear part of the Navier-Stokes equation; so, in this proof, we will be analysing the term coming from the stochastic perturbation and, for simplicity, we omit the index $s$ in the corresponding Christoffel symbols.

Writing
\begin{align*}
\nabla^s_{e_k} \nabla^s_{e_k} \tilde{u}  & = \partial_{e_k} (\nabla^s_{e_k} \tilde{u})+\Gamma_{k,l} (\nabla^s_{e_k} \tilde{u})^l \\
&=\partial_{e_k}\partial_{e_k}\tilde{u} +\Gamma_{k,l} \Gamma_{k,m}^{l} \tilde{u}^{m}
\end{align*}
we have
\begin{align*}
&\sum_{0<|k|\leq N}\nabla^s_{e_k} \nabla^s_{e_k} \tilde{u} -\partial_{e_k}\partial_{e_k} \tilde{u}-R^{N}_{s} (\tilde{u})\\
&=\sum_{0<|k|\leq N}\sum_{l,m}\Gamma_{k,l}[e_k ,e_m ]^l \tilde{u}^m -\Gamma_{l,k}[e_k ,e_m ]^l \tilde{u}^m\\
&=\sum_{0<|k|\leq N} \sum_{l,m}[e_k ,e_l ] [e_k ,e_m ]^l \tilde{u}^m
\end{align*}

\indent We analyse the component in $B^s_m$ (the one in $A^s_m$ will give a similar result).

\indent We have to sum the last expression in the indices $k$ and $l$. The two non zero contributions come from
\begin{align*}
I&= [A^s_k ,B^s_m ]^{A^s_l} [A^s_k , A^s_l ] \tilde{u}^{B^s_m}\\
II&= [B^s_k ,B^s_m ]^{B^s_l} [B^s_k , B^s_l ] \tilde{u}^{B^s_m}
\end{align*}

\indent Computation of $I$:

\indent When $l=k+m$,
\begin{align*}
I=& -\frac{[k,m]}{2 |k|^{s+1} |m|^{s+1}}|k+m|^{s+1} \frac{[k,k+m]}{2|k|^{s+1} |k+m|^{s+1}} (|2k+m|^{s+1} B^s_{2k+m}\\
& +|m|^{s+1} B^s_m )\tilde{u}^{B^s_m}\\
=& -\frac{[k,m]^2}{4 |k|^{2(s+1)} |m|^s}(|2k+m|^{s+1} B^s_{2k+m} +|m|^{s+1} B^s_m )\tilde{u}^{B^s_m}
\end{align*}
and, when $l=k-m$,
\begin{align*}
I=& -\frac{[k,m]}{2 |k|^{s+1} |m|^{s+1}}|k-m|^{s+1} \frac{[k,m]}{2|k|^{s+1} |k-m|^{s+1}} (|2k-m|^{s+1} B^s_{2k-m} \\
&+|m|^{s+1} B^s_m )\tilde{u}^{B^s_m}\\
=& -\frac{[k,m]^2}{4 |k|^{2(s+1)} |m|^{s+1}}(|2k-m|^{s+1} B^s_{2k-m} +|m|^{s+1} B^s_m )\tilde{u}^{B^s_m}
\end{align*}
(notice that $B^s_{-m}=B^s_m$).

\indent Concerning $II$, we have, for $l=k+m$,
\begin{equation*}
II= \frac{[k,m]^2}{4|k|^{2(s+1)} |m|^{s+1}} (|2k+m|^{s+1} B^s_{2k+m} -|m|^{s+1} B^s_m )\tilde{u}^{B^s_m}
\end{equation*}
and, for $l=k-m$,
\begin{equation*}
II= \frac{[k,m]^2}{4|k|^{2(s+1)} |m|^{s+1}} (|2k-m|^{s+1} B^s_{2k-m} -|m|^s B^s_m )\tilde{u}^{B^s_m}
\end{equation*}

Summing up in $k$ the non diagonal terms cancel and 
\begin{align*}
&\sum_{0<|k|\leq N} (\nabla^s_{e_k} \nabla^s_{e_k}\tilde{u} )^{B^s_m}- (R^N_s \tilde{u} )^{B^s_m}
=-\sum_k \frac{[k,m]^2}{2|k|^{2(s+1)}} \tilde{u}^{B^s_m}\\
&=-\sum_k \frac{k_1^2 m_2^2}{2|k|^{2(s+1)}}\tilde{u}^{B^s_m}-\sum_k \frac{k_2^2 m_1^2}{2|k|^{2(s+1)}}\tilde{u}^{B^s_m}+\sum_k \frac{k_1 m_2 k_2 m_1}{2|k|^{2(s+1)}} \tilde{u}^{B^s_m}
\end{align*}

Since $c_N=\sum_{0<|k|\leq N} \frac{k_1^2}{|k|^{2(s+1)}}=\sum_{0<|k|\leq N} \frac{k_2^2}{|k|^{2(s+1)}}$,
\begin{equation*}
\sum_{0<|k|\leq N} \frac{1}{2}((\nabla^s_{e_k} \nabla^s_{e_k}\tilde{u} )^{B^s_m}-\frac{\nu}{c_N}(R^N_s \tilde{u})^{B^s_m})=-\nu  m^2 \tilde{u}^{B^s_m}
\end{equation*}
that is, we have the spectral decomposition of the Laplace operator, and the result follows. 

\indent The first statement of the theorem, namely  the characterization of $u(g^N_u )$ as a martingale, comes from It\^o formula and the fact that the bounded variation part of this process is precisely the expression $\frac{\partial \tilde u}{\partial t} +\tilde u^k \nabla_{e_k} \tilde u - \L^{s,N} \tilde{u} +\frac{\nu}{2c_N } R^N_s (\tilde u)$.
\end{proof}

\section{Final comments}

The probabilistic characterization  of the Navier-Stokes solutions described above allows to
study various of their properties. But one can legitemely ask wether it shares any light on the
 construction of the solutions itself. This would amount to ask wether we can solve directly 
 (and in a probabilistic way) the stochastic equation
\begin{equation*}
DDg^N_u (t)+\frac{\nu}{2c_N } R^N_s (D g^N_u )=0 
\end{equation*}
Stated like this, we do not have a straightforward answer. Nevertheless one
can use another probabilistic characterization, in terms of forward-backward
stochastic equations, based on the observation that the
process $g^N_u (t)$ satisfies not only
\begin{equation*}
d g^N_u (t) =(u_t dt +\sqrt\frac{\nu}{c_N } od x_N (t))g^N_u (t), ~~~ g^N_u (0)=e
\end{equation*}
but also the equation
\begin{equation*}
d u_t =-\nabla p (t, g^N_u (t)) dt + X_t dx_N (t), ~~~ Dg^N_u (T) = u(T, g^N_u (T))
\end{equation*}

\indent Considering these two equations as a system on $(g^N_u (t), u_t , X_t )$ we have a forward-backward system to solve (analytically and numerically) (cf. \cite{C-S} and work in progress).

\medskip


\begin{thebibliography}{1}





\bibitem{A-M} H. Airault, P. Malliavin, \textit{Quasi-invariance of Brownian measures on the group of circle homeomorphisms and infinite-dimensional Riemannian geometry.}J. Funct. Anal.
\textbf{241} (2006), 99--142.


\bibitem{A-C} A. Albeverio, A.B. Cruzeiro, \textit{Global flows with invariant (Gibbs) measures for Euler and Navier-Stokes two dimensional fluids.} Comm. Math. Phys. 
\textbf{129} (1990), 431--444.

\bibitem{A} V.I. Arnold, \textit{Sur la g\'eom\'etrie diff\'erentielle des groupes de Lie de dimension infinie et ses aplications \`a l' hidrodynamique des fluides parfaits.} Ann. Inst. Fourier
\textbf{16} (1966), 316--361.

\bibitem{B} G.K. Batchelor, \textit{The theory of homogeneous turbulence.} Cambridge Monographs on Mechanics \& Applied Math. Series, 1953.



\bibitem{C-C} F. Cipriano, A.B. Cruzeiro, \textit{Navier-Stokes equation and diffusions on the group of homeomorphisms of the torus.} Comm. Math. Phys.
\textbf{275} (2007), 255--269.

\bibitem{C-F-M} A.B. Cruzeiro, F. Flandoli, P. Malliavin, \textit{Brownian motion of volume preserving diffeomorphisms and existence of global solutions of 2D stochastic Euler equation.} J. Funct. Anal.
\textbf{242} (2007), 304--326.

\bibitem{C-M} A.B. Cruzeiro, P. Malliavin, \textit{Nonergodicity of Euler fluid dynamics on tori versus positivity of the Arnold-Ricci tensor.} J. Funct. Anal.
\textbf{254} (2008), 1903--1925.

\bibitem{C-S} A.B. Cruzeiro, E. Shamarova, \textit{Navier-Stokes equations and forward-backward SDEs on the group of diffeomorphisms of a torus.} http://arXiv.org/abs/0807.0421.


\bibitem{E-M} D. Ebin, J. Marsden, \textit{Groups of diffeomorphisms and the motion of an incompressible fluid.} Ann.  of Math. 
\textbf{92}(1) (1970), 102--163.

\bibitem{F} A.V. Fursikov, A.I. Komech, M.I. Vishik, \textit{Some mathematical problems of statistical hydrodynamics.} Dokl. Akad. Nauk. SSSR
\textbf{246}(5) (1979), 1037--1041.

\bibitem{G} R. Galanti, A. Tsinober, \textit{Is turbulence ergodic?.} Phys. Letters A \textbf{330} (2004), 173--180.

\bibitem{K} H. Kunita, \textit{Stochastic flows and stochastic differential equation.} Cambridge University Press, 1990.

\bibitem{M} P. Malliavin, \textit{The canonic diffusion above the diffeomorphism group of the circle.} C.R. Acad. Sci. Paris \textbf{329} Ser. I (1999), 325--329.

\bibitem{N} T. Nakagomi, K. Yasue, J.C. Zambrini, \textit{Stochastic variational derivations of the Navier-Stokes equation.} Lett. Math. Phys. \textbf{160} (1999), 337--365.

\bibitem{R} R. Robert, \textit{L'effet papillon n'existe plus!.} SMF-Gazette-90, 2001.

\bibitem{S} S. Skholler, \textit{Geometry and curvature of diffeomorphism group with $H^1$ metric and mean hydrodynamics.} J. Funct. Anal. \textbf{160}(1) (1998), 337--365.

\bibitem{Y} K. Yasue, \textit{A variational principle for the Navier-Stokes equation.} J. Funct. Anal. \textbf{51}(2) (1983), 133--141.

\end{thebibliography}
\end{document}